%% file: IEEE_general.tex
\documentclass[10pt,conference,letterpaper]{IEEEtran}

\usepackage{graphicx}
\usepackage{subfigure} 
\usepackage{algorithm} % For algorithms
\usepackage{algorithmic}

\usepackage{amsmath}
\usepackage{amssymb}
\usepackage{amsthm}

\newcommand{\vw}{ {\vect{w}} }
\newcommand{\vg}{ {\vect{g}} }
\newcommand{\step}{ {\alpha} }
\newcommand{\poly}{ \text{P} }
\newcommand{\gen}{ \text{G} }
\newcommand{\Lgen}{ {\eta} }

\newtheorem{assumption}{Assumption}
\newtheorem{lemma}{Lemma}
\newtheorem{theorem}{Theorem}

\input{SuchaLib.tex}

\begin{document} 
\title{Staggered Time Average Algorithm for Stochastic Non-smooth Optimization with $O(1/T)$ Convergence}

\author{
  \IEEEauthorblockN{Sucha Supittayapornpong,~~Michael J. Neely}
  \IEEEauthorblockA{Department of Electrical Engineering\\
    University of Southern California\\
    Email: supittay@usc.edu,~~mjneely@usc.edu}
}

\maketitle

\begin{abstract} 
Stochastic non-smooth convex optimization constitutes a class of problems in machine learning and operations research.  This paper considers minimization of a non-smooth function based on stochastic subgradients.  When the function has a \emph{locally polyhedral} structure, a staggered time average algorithm is proven to have $O(1/T)$ convergence rate.  A more general convergence result is proven when the locally polyhedral assumption is removed.  In that case, the convergence bound depends on the curvature of the function near the minimum.  Finally, the locally polyhedral assumption is shown to improve convergence beyond $O(1/T)$ for a special case of deterministic problems.
\end{abstract}

\section{Introduction}
Non-smooth convex optimization constitutes a class of problems in machine learning and operations research.  Examples include optimization of the hinge loss function \cite{Lorenzo:Hinge} (used in support vector machines) and $1$-norm regularization in regression problems \cite{Tibshirani:Lasso}.  A non-smooth function is continuous and non-differentiable \cite{Nesterov:Intro}.  This lack of differentiability makes it challenging to design algorithms with fast convergence.

This paper considers a stochastic optimization problem:
\begin{equation}
  \label{problem:main}
  \min_{\vw \in \set{W}} F(\vw)
\end{equation}
where $\set{W}$ is a closed and convex set and $F: \set{W} \rightarrow \RealSet$ is a continuous, convex (but not necessary strongly convex), and non-smooth function.  Function $F$ may not be known.  The optimization proceeds by obtaining an unbiased stochastic subgradient of $F$ from an oracle.  This model with an oracle has been previously used in literature, such as \cite{Rakhlin:averaging,Ohad:averaging}.  Formally, let $\vg(\vw)$ be the subgradient of $F$ at $\vw$.  Receiving $\vw \in \set{W}$, the oracle gives an unbiased stochastic subgradient $\hat{\vg}(\vw)$ of $F$ at $\vw$ satisfying $\expect{ \hat{\vg}(\vw) | \vw } = \vg(\vw)$.  Note that $\vg(\vw)$ may not be known, but $\hat{\vg}(\vw)$ is known.  An algorithm proceeds by generating a sequence of $\vw_t$ vectors that are given to the oracle.  The next vector $\vw_{t+1}$ is determined as a function of the history of oracle outputs.  The history and the $\prtc{ \vw_t }$ sequence is used to compute an estimate $\hat{\vw}$ of the optimal solution.  It has been noted in \cite{Ohad:averaging} that this model can be applied to a class of learning problems in \cite{Shai:stocvopt}.%, where $\set{W}$ is a hypothesis class and a set of independent and identically distributed samples, that one wants to fine a predictor $\hat{\vw}$ whose expected loss $F(\vw)$ is close to optimal over $\set{W}$.

Given $\epsilon > 0$, an estimate $\hat{\vw}$ is an $O(\epsilon)$-approximation if
\begin{equation*}
  \expect{F(\hat{\vw})} - \min_{\vw \in \set{W}} F(\vw) \leq O(\epsilon).
\end{equation*}
Let $T$ be the number of unbiased stochastic subgradients obtained from an oracle.  The convergence rate is determined by the rate at which the estimate converges to the true answer, as a function of $T$.  For example, an algorithm with $O(1/\sqrt{T})$ convergence rate provides the estimate whose deviation from optimality decays to zero like:
\begin{equation*}
  \expect{F(\hat{\vw})} - \min_{\vw \in \set{W}} F(\vw) \leq O(1/\sqrt{T}).
\end{equation*}

Prior work in \cite{Boyd:stograd,Ohad:averaging} develops algorithms with $O(1/\sqrt{T})$ convergence rate when the function is non-smooth.  A smoothing method in \cite{Hua:smoothing} improves the convergence rate to $O(1/T)$ when $F$ is a linear combination of a smooth convex function and a non-smooth convex function with special structure.  Related improvements can be shown when the non-smooth function is \emph{strongly convex} \cite{Rakhlin:averaging,Ohad:averaging,Hazan:Beyond}.  The suffix algorithm in \cite{Rakhlin:averaging} is shown to have $O(1/T)$ convergence rate.  The work by \cite{Ohad:averaging} shows that using the reducing step size leads to $O(\log(T)/T)$ and $O(1/T)$ convergence rates for the last round solution and a solution calculated from a polynomial-decay averaging.  Then, the algorithm achieving optimal $O(1/T)$ convergence rate is developed in \cite{Hazan:Beyond}.  Note that all previous results that achieve the $O(1/T)$ convergence rate rely on either a restrictive strong convexity or a special structure of $F$.  It has become an open problem whether the $O(1/T)$ convergence rate can be achieved for a non-smooth and non-strongly convex function.

In this paper, a non-smooth convex function is considered.  In the case when the function has a \emph{locally polyhedral} structure \footnote{The locally polyhedral structure is also called weak sharp minima in previous literature \cite{Burke:WSM}.  It is also the generalization of a sharp minimum function in \cite{Polyak:Intro}}, a staggered time average algorithm, based on a stochastic subgradient algorithm with constant step size, is proposed.  The algorithm calculates the $O(\epsilon)$-approximation estimates with $O(1/T)$ convergence rate.  For a general convex function, the convergence rate depends on the curvature near the minimum.  To our knowledge, with the locally polyhedral structure, this is the first $O(1/T)$ convergence rate for a non-smooth convex function, which requires neither strong convexity nor smoothing.

%Exploiting these structures,   These phases are studied in their own theoretical interest before designing mechanisms to obtain the $O(\epsilon)$-approximation solution.  When the function satisfies \emph{locally-polyhedral} structure, the expected transient time is proven to be $O(1/\step)$.  When the function satisfies a more general \emph{locally-non-polyhedral} structure with parameters $S$ and $L(S)$, the expected transient time becomes $O(1/(\step L(S)))$.  Then estimated solutions are calculated from the averages of solutions obtained in the steady state phase.  They converge as $O(\step/T + \step)$ and $O(\phi(\step, S)/T + \step)$ respectively under the locally-polyhedral and locally-non-polyhedral structures, where $\phi(\step, S)$ is derived from the latter structure.  

%In retrospective, these convergence results exploit the structures of functions to choose appropriate constant step size instead of the greedily reducing step size.  A discussion on the choice of the step size to achieve the $O(\epsilon)$-approximation estimate is provided in the paper.

The paper is organized as follows.  Section \ref{sec:pre} provides notations, the staggered time average algorithm, and preliminary results.  Sections \ref{sec:poly} and \ref{sec:gen} prove respectively the results under the locally polyhedral and the general convex structures.  Section \ref{sec:deterministic} shows a fast convergence for deterministic problems.  Experiments are performed in Section \ref{sec:sim}.  Section \ref{sec:conclusion} concludes the paper.

\section{Preliminaries}
\label{sec:pre}
The closed convex set $\set{W}$ is a subset of $\RealSet^N$, for some positive integer $N$, with Euclidean norm $\norm{\cdot}$ and inner product $\inner{\cdot}{\cdot}$.  Function $F$ is assumed to be convex (possibly non-smooth) over $\set{W}$ and satisfies the following assumption.  Define $F^\ast \defeq \inf_{\vw \in \set{W}} F(\vw)$.
\begin{assumption}
  \label{ass:solution}
  The minimum of $F$ is achievable in $\set{W}$, and the set of optimal solutions $\set{W}^\ast = \Arginf_{\vw \in \set{W}} F(\vw) = \prtc{ \vw^\ast \in \set{W} : F(\vw^\ast) = F^\ast }$ is closed.
\end{assumption}

A subgradient of $F$ at $\vw \in \set{W}$ is denoted by $\vg(\vw)$ and satisfies for any $\vw' \in \set{W}$:
\begin{equation}
  \label{eq:Fconvex}
  F(\vw') \geq F(\vw) + \inner{ \vg(\vw) }{\vw' - \vw}.
\end{equation}
An unbiased stochastic subgradient at $\vw$ is denoted by $\hat{\vg}(\vw)$, which satifies $\expect{\hat{\vg}(\vw)| \vw} = \vg(\vw)$.
\begin{assumption}
  \label{ass:G}
There exists a constant $G < \infty$ such that
\begin{equation*}
  \norm{\hat{\vg}(\vw)} \leq G \quad\quad \forall \vw \in \set{W}.
\end{equation*}
\end{assumption}
Assumption \ref{ass:G} is also used in previous literature.

A stochastic subgradient algorithm with a positive constant step size $\step > 0$ initializes $\vw_0 \in \set{W}$ and proceeds repeatedly as
\begin{equation}
  \label{alg:subgrad}
  \vw_{t+1} = \Pi_{\set{W}} \prts{ \vw_t - \step \hat{\vg}(\vw_t) } \quad \forall t \in \prtc{0, 1, 2, \dotsc},
\end{equation}
where $\Pi_\set{W}$ denotes projection on $\set{W}$ and $\vw_t$ denotes the values of $\vw$ at round $t$.

\subsection{Staggered Time Averages}
The staggered time average algorithm is summarized in Algorithm \ref{alg:main}.  Define:
\begin{equation*}
  \bar{\vw}_{2^k-1}^T = \frac{1}{T} \sum_{t = 0}^{T-1} \vw_{2^k-1+t} .
\end{equation*}
This average can be computed on-the-fly as shown in Algorithm \ref{alg:main}.

\begin{algorithm}
\caption{Staggered Time Averages}
\label{alg:main}
\begin{algorithmic}
  \STATE Initialize: $\vw_0 \in \set{W},~ \step > 0$
  \FOR{$ t \in \prtc{0, 1, 2, \dotsc}$}
    \STATE {\it // Staggered time averages}
    \IF{$t = 2^{k^\ast} - 1$ for some $k^\ast \in \prtc{0, 1, 2, \dotsc}$}
      \STATE $k \leftarrow k^\ast$
      \STATE $\bar{\vw}_{2^k - 1}^1 \leftarrow \vw_t$
    \ELSE
      \STATE $\bar{\vw}_{2^k - 1}^{t-2^k+2} \leftarrow \frac{t-2^k+1}{t-2^k+2}\bar{\vw}_{2^k - 1}^{t-2^k+1} + \frac{1}{t-2^k+2}\vw_t$
    \ENDIF
    \STATE {\it // Stochastic subgradient}
    \STATE $\vw_{t+1} \leftarrow \Pi_{\set{W}}\prts{ \vw_t - \step \hat{\vg}(\vw_t) }$
  \ENDFOR
\end{algorithmic}
\end{algorithm}

Algorithm \ref{alg:main} implements the subgradient algorithm \eqref{alg:subgrad} with constant step size in each round.  The staggered time averages reset the calculation of estimates every $2^k-1$ for $k \in \prtc{0, 1, 2, \dotsc}$.  Specifically, for every $k \in \prtc{0, 1, 2, \dotsc}$, the algorithm generates estimates $\bar{\vw}_{2^k-1}^T$ for $T \in \prtc{1, \dotsc, 2^k}$.  To analyze Algorithm \ref{alg:main}, the properties of the subgradient algorithm are proven in this section.  Then the staggered time averages are analyzed in Section \ref{sec:poly} and Section \ref{sec:gen}.

Note that Algorithm \ref{alg:main} is different from the suffix averaging in \cite{Rakhlin:averaging} which uses the reducing step size.

\subsection{Basic Results}
We consider algorithm \eqref{alg:subgrad} with a positive constant step size $\step$.  The initial $\vw_0 \in \set{W}$ is any constant vector.  For every $\vw_t \in \set{W}$, define the closest optimal solution to $\vw_t$ as
\begin{equation*}
  \vw_t^\ast = \arginf_{\vw \in \set{W}^\ast} \norm{ \vw - \vw_t }.
\end{equation*}
Under Assumption \ref{ass:solution}, this $\vw_t^\ast$ is unique because of the convexity and the closeness of $\set{W}^\ast$.  The following lemma modifies a well known manipulation.
\begin{lemma}
  \label{lem:basic}
Suppose Assumptions \ref{ass:solution} and \ref{ass:G} hold.  It holds for any $t \in \prtc{0, 1, 2, \dotsc}$ that
\begin{multline*}
  \expect{ \norm{\vw_{t+1} - \vw_{t+1}^\ast}^2 \Big| \vw_t } \leq \norm{\vw_t - \vw_t^\ast}^2 + \step^2G^2 \\+ 2\step\prts{ F^\ast - F(\vw_t) }.
\end{multline*}
\begin{proof}
For any $t$, by definition of $\vw_{t+1}^\ast$ as the minimizer of $\norm{\vw - \vw_{t+1}}^2$ over all $\vw \in \set{W}$, we have:
\begin{align*}
  \norm{\vw_{t+1} - \vw_{t+1}^\ast}^2
  & \leq \norm{ \vw_{t+1} - \vw_t^\ast }^2 \\
  & = \norm{ \Pi_{\set{W}}\prts{ \vw_t - \step\hat{\vg}(\vw_t) } - \vw_t^\ast } \\
  & \leq \norm{ \vw_t - \step\hat{\vg}(\vw_t) - \vw_t^\ast }^2
\end{align*}
where the final inequality holds by the non-expansive property of the projection onto the convex set $\set{W}$ \cite{Bertsekas:Convex}.  Expanding the right-hand-side gives:
\begin{align*}
  & \norm{\vw_{t+1} - \vw_{t+1}^\ast}^2  \\
  & \quad = \norm{ \vw_t - \vw_t^\ast }^2 + \step^2\norm{ \hat{\vg}(\vw_t) }^2 - 2\step\inner{ \hat{\vg}(\vw_t) }{ \vw_t - \vw_t^\ast } \\
  & \quad \leq \norm{ \vw_t - \vw_t^\ast }^2 + \step^2G^2 + 2\step\inner{ \hat{\vg}(\vw_t) }{ \vw_t^\ast - \vw_t }
\end{align*}

Taking a conditional expectation given $\vw_t$ yields
\begin{multline*}
  \expect{\norm{ \vw_{t+1} - \vw_{t+1}^\ast }^2 \Big| \vw_t} \leq \norm{ \vw_t - \vw_t^\ast }^2 + \step^2G^2 \\+ 2\step\inner{ \vg(\vw_t) }{ \vw_t^\ast - \vw_t }.
\end{multline*}
Using the subgradient property in \eqref{eq:Fconvex} and the fact that $F(\vw_t^\ast) = F^\ast$ proves the lemma.
\end{proof}
\end{lemma}
Note that Lemma \ref{lem:basic} uses a projection technique similar to standard analysis for the subgradient projection algorithm, (as in \cite{Zinkevich:OnlineConvex}, \cite{Nesterov:Intro} or \cite{Bertsekas:Convex}).  The standard approach compares the current iterate to a fixed optimal point $\vw^\ast$.  Lemma \ref{lem:basic} compares to the closest point in the optimal set.  This is a simple but important distinction that is crucial in later sections for improved convergence time results.

While Lemma \ref{lem:basic} is stated in a form useful for the analysis of later sections, it can readily be used to establish the standard $O(1/\epsilon^2)$ result for convex functions (see also \cite{Zinkevich:OnlineConvex,Nesterov:Intro,Bertsekas:Convex}).

Define an average of $T$-consecutive solutions from $t_0$ as
\begin{equation}
  \label{eq:avg_w}
  \bar{\vw}_{t_0}^T \defeq \frac{1}{T} \sum_{t = t_0}^{t_0 + T -1} \vw_t.
\end{equation}

Taking an expectation of the result in Lemma \ref{lem:basic} gives
\begin{multline*}
  \expect{ F(\vw_t) } - F^\ast \leq \\
  \frac{\step G^2}{2} + \frac{1}{2\step}\expect{ \norm{ \vw_{t} - \vw_t^\ast }^2 - \norm{ \vw_{t+1} - \vw_{t+1}^\ast }^2 }.
\end{multline*}
Summing from $t_0$ to $t_0+T-1$ and dividing by $T$ gives
\begin{multline*}
  \frac{1}{T} \sum_{t=t_0}^{t_0+T-1} \expect{ F(\vw_t) } - F^\ast \leq \\ \frac{\step G^2}{2} + \frac{1}{2\step T}\expect{ \norm{ \vw_{t_0} - \vw_{t_0}^\ast }^2 - \norm{ \vw_{t_0+T} - \vw_{t_0+T}^\ast }^2 }.
\end{multline*}
Using Jensen's inequality and convexity of $F$, definition \eqref{eq:avg_w}, and non-negativity of $\norm{ \vw_{t_0+T} - \vw_{t_0+T}^\ast }^2$ yield:
\begin{equation}
  \label{eq:avgbasic}
  \expect{ F(\bar{\vw}_{t_0}^T) } - F^\ast \leq \frac{\step G^2}{2} + \frac{1}{2\step T} \expect{ \norm{ \vw_{t_0} - \vw_{t_0}^\ast }^2 },
\end{equation}
for any $t_0 \in \prtc{0, 1, 2, \dotsc}$ and any positive integer $T$.

Equation \eqref{eq:avgbasic} suggests that, to achieve an $O(\epsilon)$-approximation estimate, one can choose step size $\step = \Theta(\epsilon)$, number of rounds $T = \Theta(1/\epsilon^2)$, and define $\hat{\vw}$ as the average of $\vw_t$ values over the first $T$ rounds.  This is equivalent to $O(1/\sqrt{T})$ convergence rate.  Fortunately, this convergence rate can be improved by starting the average at an appropriate time depending on the structure of the function $F$.  These structures are shown in Section \ref{sec:poly} and Section \ref{sec:gen}.  We first prove several useful results used in those sections.

\subsection{Concentration Bound}
These results are used to upper bound $\expect{\norm{ \vw_{t_0} - \vw_{t_0}^\ast }^2}$ in \eqref{eq:avgbasic}.  Define $K_t \defeq \norm{ \vw_t - \vw_t^\ast }$ for all $t \in \prtc{0, 1, 2, \dotsc}$.
\begin{lemma}
  \label{lem:k}
Suppose Assumption \ref{ass:G} holds.  It holds for every $t \in \prtc{0, 1, 2, \dotsc}$ that
\begin{align}
  \abs{ K_{t+1} - K_t } & \leq 2 \step G \label{eq:abs_bound}\\
  \expect{ K_{t+1} - K_t | \vw_t } & \leq 2 \step G. \label{eq:drift_bound}
\end{align}
\begin{proof}
The first part is proven in two cases.\\
i) If $K_{t+1} \geq K_t$, definition of $\vw_{t+1}$ in \eqref{alg:subgrad} and the non-expansive projection implies
\begin{align*}
  & \abs{ K_{t+1} - K_t } = K_{t+1} - K_t \leq \norm{ \vw_{t+1} - \vw_{t}^\ast } - K_t \\
  & \quad \leq \norm{ \vw_{t} - \step \hat{\vg}(\vw_t) - \vw_t^\ast } - K_t \\
  & \quad \leq \norm{ \vw_{t} - \vw_t^\ast } + \step \norm{\hat{\vg}(\vw_t)} - K_t \leq \step G.
\end{align*}
ii) If $K_{t+1} < K_t$, we have
\begin{align*}
  & \abs{ K_{t+1} - K_t } = K_t - K_{t+1} \\
  & \quad = \norm{ \vw_t - \vw_{t+1} + \vw_{t+1} - \vw_{t+1}^\ast + \vw_{t+1}^\ast - \vw_t^\ast } - K_{t+1}\\
  & \quad \leq \norm{\vw_t - \vw_{t+1}} + K_{t+1} + \norm{\vw_{t+1}^\ast - \vw_t^\ast} - K_{t+1} \\
  & \quad \leq 2 \norm{\vw_t - \vw_{t+1}} \leq 2 \step \norm{\hat{\vg}(\vw_t)} \leq 2 \step G,
\end{align*}
where the last line uses non-expansive projection implying $\norm{\vw_{t+1}^\ast - \vw_t^\ast} \leq \norm{\vw_{t+1} - \vw_t}$.

These two cases prove the first part.  The second part follows by taking a conditional expectation given $\vw_t$ of $K_{t+1} - K_t \leq \abs{K_{t+1} - K_t}$.
\end{proof}
\end{lemma}

The concentration bound reinterprets the lemma in \cite{Neely:Convergence}.
\begin{lemma}
  \label{lem:concentration}
Suppose Assumption \ref{ass:G} holds and there exists a positive real-valued $\beta$ and a $\gamma \in \RealSet$ such that for all $t \in \prtc{0, 1, 2, \dotsc}$:
\begin{equation}
  \label{eq:cond_concentration}
  \expect{K_{t+1} - K_t | K_t } \leq \left\{
  \begin{array}{rl}
    2 \step G &,\text{if}~ K_t < \gamma\\
    -\beta &,\text{if}~ K_t \geq \gamma.
  \end{array} \right.
\end{equation}
Assume $K_{t_0} = k_{t_0}$ (with probability 1) for some $k_{t_0} \in \RealSet$.  Then following holds for all $t \in \prtc{t_0, t_0+1, t_0 +2, \dotsc}$
\begin{equation}
  \label{eq:concentration}
  \expect{ e^{rK_t} } \leq D + \prtr{ e^{rk_{t_0}} - D } \rho^{t - t_0},
\end{equation}
where constants $r, \rho,$ and $D$ are:
\begin{align*}
  r & \defeq \frac{3 \beta}{12 \step^2 G^2 + 2 \step G \beta}, \quad\quad\rho \defeq 1 - \frac{ r\beta }{2}, \\
  D & \defeq \frac{ \prtr{ e^{2 \step G r} - \rho } e^{r\gamma} }{1 - \rho}.
\end{align*}
It can be shown that \eqref{eq:abs_bound} and \eqref{eq:cond_concentration} together imply $0 < \beta \leq 2 \step G$, and hence it can be shown that $0 < \rho < 1$.
\end{lemma}

Constants $\beta$ and $\gamma$ in Lemma \ref{lem:concentration} depend on a structure of function $F$.  We then look at the first structure.

\section{Locally Polyhedral Structure}
\label{sec:poly}
\begin{figure}
  \centering
  \includegraphics[scale=0.9]{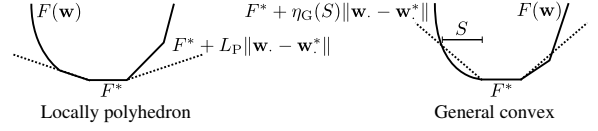}
  \caption{Structures of function $F$}
  \label{fig:structure}
\end{figure}

In this section, function $F$ is assumed to have a \emph{locally polyhedral} structure, which is illustrated in Figure \ref{fig:structure}.  This structure is generalized from \cite{Longbo:Lagrange}.  It is assumed throughout that Assumptions \ref{ass:solution} and \ref{ass:G} still hold.  Note that, in machine learning, this $F$ can be the hinge loss function with $1$-norm regularization.

\begin{assumption}
  \label{ass:poly}
(Locally polyhedral assumption) There exists a constant $L_\poly > 0$ such that for every $\vw_t \in \set{W}$ the following holds
\begin{equation}
  \label{eq:structure_poly}
  F(\vw_t) - F^\ast \geq L_\poly \norm{\vw_t - \vw_t^\ast}.
\end{equation}
\end{assumption}
The subscript ``P'' in $L_\poly$ represents ``polyhedral.''

\subsection{Drift and Transient Time}
Using this locally polyhedral structure, sequence $\prtc{ \vw_t }_{t=0}^\infty$ generated by algorithm \eqref{alg:subgrad} has the following drift property.  Define 
\begin{equation*}
  B_\poly \defeq \max\prts{ \frac{L_\poly}{2}, \frac{G^2}{L_\poly} }.
\end{equation*}

\begin{lemma}
  \label{lem:drift_poly}
  Suppose Assumptions \ref{ass:solution}, \ref{ass:G}, \ref{ass:poly} hold.  For any $\vw_t \in \set{W}$ and $\norm{ \vw_t - \vw_t^\ast } \geq \step B_\poly$, the following holds
\begin{equation}
  \label{eq:drift_poly}
  \expect{ \norm{ \vw_{t+1} - \vw_{t+1}^\ast } | \vw_t } \leq \norm{ \vw_t - \vw_t^\ast } - \frac{\step L_\poly}{2}.
\end{equation}
\begin{proof}
For any $\vw_t \in \set{W}$, if condition
\begin{equation}
  \label{eq:lempoly-1}
  \step^2 G^2 + 2\step \prts{ F^\ast - F(\vw_t) } \leq -\step L_\poly \norm{ \vw_t - \vw_t^\ast } + \frac{\step^2 L_\poly^2}{4}
\end{equation}
is true, then the result in Lemma \ref{lem:basic} implies that
\begin{align*}
  & \expect{ \norm{ \vw_{t+1} - \vw_{t+1}^\ast }^2 \Big| \vw_t } \\
  & \quad \leq \norm{ \vw_t - \vw_t^\ast }^2 - \step L_\poly \norm{ \vw_t - \vw_t^\ast } + \frac{\step^2 L_\poly^2}{4} \\
  & \quad = \prtr{ \norm{ \vw_t - \vw_t^\ast } - \frac{\step L_\poly}{2} }^2.
\end{align*}
Applying Jensen's inequality on the left-hand-side gives
\begin{equation*}
  \prtr{ \expect{ \norm{ \vw_{t+1} - \vw_{t+1}^\ast } | \vw_t } }^2 \leq \prtr{ \norm{ \vw_t - \vw_t^\ast } - \frac{\step L_\poly}{2} }^2.
\end{equation*}
When $\norm{ \vw_t - \vw_t^\ast } \geq \step B_\poly \geq \step L_\poly/2$, inequality \eqref{eq:drift_poly} holds.  

It remains to show that condition \eqref{eq:lempoly-1} must hold whenever $\norm{ \vw_t - \vw_t^\ast } \geq \step B_\poly$.  Starting with the left-hand-side of \eqref{eq:lempoly-1} we have, for every $\vw_t \in \set{W}$:
\begin{equation*}
  \step^2 G^2 + 2\step \prts{ F^\ast - F(\vw_t) } \leq \step^2 G^2 - 2 \step L_\poly \norm{ \vw_t - \vw_t^\ast }
\end{equation*}
where the inequality follows the locally polyhedral structure \eqref{eq:structure_poly}.  By the definition of $B_\poly$, it holds that $\norm{ \vw_t - \vw_t^\ast } \geq \step G^2/L_\poly$ and $\vw_t \in \set{W}$.  Substituting this into the above gives the result.
\end{proof}
\end{lemma}

This lemma implies that, when the distance between $\vw_t$ and $\vw_t^\ast$ is at least $\step B_\poly$, then $\vw_{t+1}$ is expected to get closer to $\vw_{t+1}^\ast$.  This phenomenon suggests that $\vw_t$ will concentrate around $\set{W}^\ast$ after some transient time, if it is not inside the set already.

%Combining Lemma \ref{lem:drift_poly} and Lemma \ref{lem:concentration} gives the transient time and a useful bound. 
Let constants $U_\poly, r_\poly, \rho_\poly, D_\poly$ be
\begin{align}
  U_\poly &\defeq \frac{2 (D_\poly + 1) }{r_\poly^2}, \label{eq:U_poly} \\
  r_\poly &\defeq \frac{3 L_\poly}{24G^2 + 2 L_\poly G}, \label{eq:r_poly} \\
  \rho_\poly &\defeq 1 - \frac{3 L_\poly^2}{4\prtr{24G^2 + 2 L_\poly G}}, \label{eq:rho_poly} \\
  D_\poly &\defeq \frac{ \prtr{ e^{2 G r_\poly} - \rho_\poly } e^{r_\poly B_\poly} }{1 - \rho_\poly}. \label{eq:D_poly}
\end{align}
Given any $\vw_0 \in \set{W}$, define
\begin{equation}
  \label{eq:T_poly}
  T_\poly(\step) \defeq \left\lceil \frac{ r_\poly \norm{ \vw_0 - \vw_0^\ast } }{ \step \log (1/ \rho_\poly ) } \right\rceil,
\end{equation}
where constants $r_\poly$ and $\rho_\poly$ are defined in \eqref{eq:r_poly} and \eqref{eq:rho_poly}.  This $T_\poly(\step)$ can be called a transient time for the locally polyhedral structure, since a useful bound holds after this time.

\begin{lemma}
  \label{lem:expect_poly}
Suppose Assumptions \ref{ass:solution}, \ref{ass:G}, \ref{ass:poly} hold.  When $t \geq T_\poly(\step)$, the following holds
\begin{equation*}
  \expect{ \norm{ \vw_t - \vw_t^\ast }^2 } \leq \step^2 U_\poly,
\end{equation*}
where constant $U_\poly$ is defined in \eqref{eq:U_poly}.
\begin{proof}
Lemma \ref{lem:drift_poly} gives $\expect{ K_{t+1} - K_t | K_t } \leq -\step L_\poly/2$ as in \eqref{eq:drift_poly} when $K_t = \norm{\vw_t - \vw_t^\ast} \geq \step B_\poly$.  Therefore, the constants $\beta$ and $\gamma$ in Lemma \ref{lem:concentration} can be set as $\beta = \step L_\poly/2$ and $\gamma = \step B_\poly$.  When $t_0 = 0$, we have $k_{t_0} = K_0 = \norm{ \vw_0 - \vw_0^\ast }$ (with probability 1).  From \eqref{eq:concentration}, it holds for all $t \in \prtc{0, 1, 2, \dotsc}$ that
\begin{align}
  \expect{ e^{\frac{r_\poly K_t}{\step}} } &\leq D_\poly + \prtr{ e^{\frac{r_\poly k_{t_0}}{\step} } - D_\poly } \rho_\poly^\prtr{ t - 0 } \notag \\
  &\leq D_\poly + e^{\frac{r_\poly K_0}{\step}}\rho_\poly^t, \label{lem-bound-2}
\end{align}
where constants $r_\poly, \rho_\poly, D_\poly$ are defined in \eqref{eq:r_poly}--\eqref{eq:D_poly} respectively.  We then show that 
\begin{equation}
  \label{lem-bound-1}
  e^{\frac{r_\poly K_0}{\step}}\rho_\poly^t \leq 1 \quad\quad \forall t \geq T_\poly(\step).
\end{equation}
Inequality $e^{\frac{r_\poly K_0}{\step}}\rho_\poly^t \leq 1$ is equivalent to $t \geq \frac{ r_\poly K_0 }{ \step \log(1/ \rho_\poly ) }$ by arithmetic and the fact that $\log(1/ \rho_\poly) > 0$.  From the definition of $T_\poly(\step)$ in \eqref{eq:T_poly}, it holds that $T_\poly(\step) \geq \frac{ r_\poly K_0 }{ \step \log(1/ \rho_\poly ) }$, and the results \eqref{lem-bound-1} follows.

From \eqref{lem-bound-1}, inequality \eqref{lem-bound-2} becomes
\begin{equation*}
  \expect{ e^{\frac{r_\poly K_t}{\step}} } \leq D_\poly + 1 \quad\quad \forall t \geq T_\poly(\step).
\end{equation*}
The Chernoff bounds (see, for example, \cite{Ross:Stochastic}) implies for any $m > 0$ that
\begin{align*}
  \prob{ K_t \geq m } &\leq e^{-\frac{r_\poly m}{\step} } \expect{ e^{\frac{r_\poly K_t}{\step}} }\\
  &\leq e^{-\frac{r_\poly m }{\step}} \prtr{ D_\poly + 1 } \quad\quad \forall t \geq T_\poly(\step).
\end{align*}
Using $\expect{ K_t^2 } = 2 \int_0^\infty m \prob{ K_t \geq m } dm$ and the above bound proves the lemma.
\end{proof}
\end{lemma}

The definition of $T_\poly(\step)$ in \eqref{eq:T_poly} implies that the transient time is $O(1/\step)$ under the locally polyhedral structure.  Then, Lemma \ref{lem:expect_poly} implies that $\expect{ \norm{ \vw_t - \vw_t^\ast }^2 } \leq \step^2 U_\poly$ every round $t$ after the transient time.

\subsection{Convergence Rate}
We are now ready to prove the convergence rate of the staggered time averages in Algorithm \ref{alg:main} under the locally polyhedral structure.

\begin{theorem}
  \label{thm:poly}
Suppose Assumptions \ref{ass:solution}, \ref{ass:G}, \ref{ass:poly} hold.  It holds for any $t \geq T_\poly(\step)$ and any positive integer $T$ that
\begin{equation}
  \label{eq:F_poly}
  \expect{ F(\bar{\vw}_{t}^T) } - F^\ast \leq \frac{\step G^2}{2} + \frac{\step U_\poly}{2T},
\end{equation}
where constant $U_\poly$ is defined in \eqref{eq:U_poly}.
\begin{proof}
The theorem follows from \eqref{eq:avgbasic} and Lemma \ref{lem:expect_poly} that
\begin{equation*}
  \expect{ F(\bar{\vw}_{t}^T) } - F^\ast \leq \frac{\step G^2}{2} + \frac{\step^2 U_\poly}{2\step T} \quad\quad \forall t \geq T_\poly(\step).
\end{equation*}
\end{proof}
\end{theorem}

After the transient time, Theorem \ref{thm:poly} implies that estimates, as the averages, converge as $O(\step + \step/T)$.  To obtain an $O(\epsilon)$-approximation estimate, the step size must be set to $\Theta(\epsilon)$.  Recall that the transient time \eqref{eq:T_poly} is $O(1/\step)$ and that Algorithm \ref{alg:main} resets the averages at round $2^k-1$ for $k \in \prtc{0, 1, 2, \dotsc}$.  Let $\hat{k} = \arginf_{k \in \prtc{0, 1, 2, \dotsc}} 2^{\hat{k}} -1 \geq T_\poly(\step)$ be the first reset after the transient time.  The exponential increasing implies that $2^{\hat{k}}-1 \leq 2 T_\poly(\step)$.  Therefore, the total time to obtain the estimate is $O(1/\epsilon)$, and the convergence rate is $O(1/T)$.  Note that the staggered time average algorithm is proposed, because $\norm{\vw_0 - \vw_0^\ast}$ in \eqref{eq:T_poly} can not be upper bounded if $\set{W}$ is unbounded.   Also, even though the convergence rate depends only on the step size, performing the averages also helps obtaining more accurate estimates, as shown in Section \ref{sec:sim}.

%However, this improved convergence rate does not require the strong convexity of a non-smooth convex function.  It is achieved by exploiting the locally-polyhedral structure and the constant step size.  

%Also any $\bar{\vw}_{2^k}^T$ and $T \in \prtc{1, \dotsc, 2^k}$ with $k \geq \hat{k}$ is an $O(\epsilon)$-approximation estimate.

\section{General Convex Function}
\label{sec:gen}
In this section, function $F$ is allowed to be a \emph{general convex} function, possibly one that does not satisfy the locally polyhedral structure (Assumption \ref{ass:poly}).  A general convex function is illustrated in Figure \ref{fig:structure}.  It is assumed throughout that Assumptions \ref{ass:solution} and \ref{ass:G} still hold.  

Define $\set{A}(S) = \prtc{ \vw_t \in \set{W}: \norm{\vw_t - \vw_t^\ast} = S}$, and define $S_\text{max}$ as the supremum value of $S > 0$ for which $\set{A}(S)$ is nonempty ($S_\text{max}$ is possibly infinity).  Assume that $S_\text{max} > 0$.  Convexity of the set $\set{W}$ implies that $\set{A}(S)$ is nonempty for all $S \in (0, S_\text{max})$.  For each $S \in (0, S_\text{max})$ define:
\begin{equation}
  \label{eq:LS}
  \Lgen(S) = \inf_{\vw_t \in \set{W}: \norm{\vw_t - \vw_t^\ast} = S} \frac{ \abs{ F(\vw_t) - F^\ast } }{ \norm{ \vw_t - \vw_t^\ast } }.
\end{equation}

\begin{lemma}
  \label{lem:gen}
Suppose Assumption \ref{ass:solution} holds.  If $F$ is convex and $\set{W}$ is closed, then for all $S \in (0, S_\text{max})$: \\
i) $\Lgen(S) > 0$\\
ii) Whenever $\vw_t \in \set{W}$ and $\norm{ \vw_t - \vw_t^\ast } \geq S$, it holds that
\begin{equation}
  \label{eq:structure_gen}
  F(\vw_t) - F^\ast \geq \Lgen(S) \norm{ \vw_t - \vw_t^\ast }.
\end{equation}
\begin{proof}
The first part is proven by contradiction.  Define $\set{A} = \prtc{\vw_t \in \set{W}: \norm{\vw_t - \vw_t^\ast} = S}$, which is nonempty.  Note that $\set{A}$ is compact.  Suppose $\Lgen(S) = 0$.  Function $\abs{ F(\vw_t) - F^\ast }/S$ is continuous.  The infimum of a continuous function over the compact set is achieved by a point in the set.  Thus, there is a point $\vect{y} \in \set{A}$ such that $\abs{ F(\vect{y}) - F^\ast }/S = 0$.  That is $F(\vect{y}) = F^\ast$, and $\vect{y} \in \set{W}^\ast$.  Since $\vect{y} \in \set{A}$, it also satisfies $\inf_{\vect{y}^\ast \in \set{W}^\ast} \norm{ \vect{y} - \vect{y}^\ast } = S$ and $\vect{y} \notin \set{W}^\ast$, which is a contradiction.

The second part is proven using the convexity of $F$.  Let $\vect{z} \in \set{W}$ be a vector such that $\norm{\vect{z} - \vect{z}^\ast} \geq S$ where $\vect{z}^\ast = \arginf_{\hat{\vect{z}} \in \set{W}^\ast} \norm{\vect{z} - \hat{\vect{z}}}$.  We want to show that $F(\vect{z}) - F^\ast \geq \Lgen(S) \norm{ \vect{z} - \vect{z}^\ast }$.  The convexity of the set $\set{W}$ implies that the line segment between $\vect{z}$ and $\vect{z}^\ast$ is inside $\set{W}$.  The convexity of $F$ over $\set{W}$ implies that $F$ is convex when it is restricted to this line segment.

Define $\vect{y}$ as a point on this line segment such that $\norm{\vect{y} - \vect{y}^\ast} = S$ where $y^\ast = \arginf_{\hat{y} \in \set{W}^\ast} \norm{\vect{y} - \hat{\vect{y}}}$.  Then both $\vect{y} \in \set{W}$ and $\vect{y} \in \set{A}$.  The convexity of $F$ over the line segment implies that
\begin{equation*}
  \frac{ F(\vect{z}) - F^\ast }{\norm{\vect{z} - \vect{z}^\ast}} \geq \frac{F(\vect{y}) - F^\ast}{\norm{\vect{y} - \vect{y}^\ast}} \geq \Lgen(S),
\end{equation*}
where the last inequality uses \eqref{eq:LS}.
%The second part is proven using convexity of $F$.  Let $\vw^\circ \in \set{W}$ be vector such that $\norm{\vw^\circ - \vw^\ast} \geq S$.  We want to show that $F(\vw^\circ) - F(\vw^\ast) \geq \Lgen(S) \norm{ \vw^\circ - \vw^\ast }$.  The convexity of the set $\set{W}$ implies that the line segment between $\vw^\circ$ and $\vw^\ast$ is inside $\set{W}$.  Define $\vw'$ as a point on this line segment such that $\norm{\vw' - \vw^\ast} = S$.  Then $\vw' \in \set{W}$ and $\vw' \in \set{A}$.  Further, $\norm{\vw^\ast - \vw'} + \norm{\vw' - \vw^\circ} = \norm{\vw^\ast - \vw^\circ}$.  Let $a = \norm{\vw^\ast - \vw'}/\norm{\vw^\ast - \vw^\circ}$, so $1-a = \norm{\vw' - \vw^\circ}/\norm{\vw^\ast - \vw^\circ}$, and $\vw' = a \vw^\ast + (1-a) \vw^\circ$.  The convexity of $F$ implies that $a F(\vw^\ast) + (1-a) F(\vw^\circ) \geq F(\vw')$.  Rearranging terms gives the following which proves this part, where the last inequality uses \eqref{eq:LS}.
%\begin{align*}
%  \frac{ F(\vw^\circ) - F(\vw^\ast) }{\norm{\vw^\circ - \vw^\ast}}
%  &\geq \frac{F(\vw') - F(\vw^\ast)}{\norm{\vw' - \vw^\circ}} \geq \Lgen(S).
%  &\geq \frac{F(\vw') - F(\vw^\ast)}{\norm{\vw^\ast - \vw^\circ}} \geq \Lgen(S)
%\end{align*}
\end{proof}
\end{lemma}

The difference between Assumption \ref{ass:poly} and Lemma \ref{lem:gen} is that the bound \eqref{eq:structure_gen} only holds when $\norm{\vw_t - \vw_t^\ast} \geq S$.  The choice of $S$ for a particular function $F$ affects the transient time and convergence of achieving $O(\epsilon)$-approximation estimates.  This effect does not occur with Assumption \ref{ass:poly}.

\subsection{Drift and Transient Time}
Using Lemma \ref{lem:gen}, the sequence $\prtc{ \vw_t }_{t=0}^\infty$ generated by algorithm \eqref{alg:subgrad} has the following property.  Define
\begin{equation}
  \label{eq:B_gen}
  B_\gen (\step, S) \defeq \max\prts{ \frac{\step \Lgen(S)}{2}, S, \frac{\step G^2}{\Lgen(S)}}.
\end{equation}

\begin{lemma}
  \label{lem:drift_gen}
  Suppose Assumptions \ref{ass:solution} and \ref{ass:G} hold.  Function $F$ is convex with $\Lgen(S)$ defined in \eqref{eq:LS} for all $S \in (0, S_\text{max})$.   For any $\vw_t \in \set{W}$ that $\norm{ \vw_t - \vw_t^\ast } \geq B_\gen(\step, S)$, the following holds
\begin{equation}
  \label{eq:drift_gen}
  \expect{ \norm{ \vw_{t+1} - \vw_{t+1}^\ast } | \vw_t } \leq \norm{ \vw_t - \vw_t^\ast } - \frac{\step \Lgen(S)}{2}.
\end{equation}
\begin{proof}
For any $\vw_t \in \set{W}$, if condition
\begin{multline}
  \label{eq:lemgen-1}
  \step^2 G^2 + 2\step \prts{ F^\ast - F(\vw_t) } \leq -\step \Lgen(S) \norm{ \vw_t - \vw_t^\ast } \\+ \frac{\step^2 \Lgen(S)^2}{4}
\end{multline}
is true, then the result in Lemma \ref{lem:basic} implies that
\begin{align*}
  & \expect{ \norm{ \vw_{t+1} - \vw_{t+1}^\ast }^2 \Big| \vw_t } \\
  & \quad \leq \norm{ \vw_t - \vw_t^\ast }^2 - \step \Lgen(S) \norm{ \vw_t - \vw_t^\ast } + \frac{\step^2 \Lgen(S)^2}{4} \\
  & \quad = \prts{ \norm{ \vw_t - \vw_t^\ast } - \frac{\step \Lgen(S)}{2} }^2.
\end{align*}
Applying Jensen's inequality on the left-hand-side gives
\begin{equation*}
  \prtr{ \expect{ \norm{ \vw_{t+1} - \vw_{t+1}^\ast } | \vw_t } }^2 \leq \prts{ \norm{ \vw_t - \vw_t^\ast } - \frac{\step \Lgen(S)}{2} }^2.
\end{equation*}
When $\norm{ \vw_t - \vw_t^\ast } \geq B_\gen(\step, S) \geq \step \Lgen(S)/2$, inequality \eqref{eq:drift_gen} holds.

It remains to show condition \eqref{eq:lemgen-1} must hold whenever $\norm{ \vw_t - \vw_t^\ast } \geq B_\gen(\step, S)$.  Starting with the left-hand-side of \eqref{eq:lemgen-1}, we have that
\begin{multline*}
   \step^2 G^2 + 2\step \prts{ F^\ast - F(\vw_t) } \leq \\\step^2 G^2 - 2 \step \Lgen(S) \norm{ \vw_t - \vw_t^\ast }
\end{multline*}
where the inequality follows \eqref{eq:structure_gen}, as $B_\gen(\step, S) \geq S$.  By the definition of $B_\gen(\step, S)$, it holds that $\norm{ \vw_t - \vw_t^\ast } \geq \step G^2/\Lgen(S)$.  Substituting this into the above inequality gives the result.
\end{proof}
\end{lemma}

This result is similar to Lemma \ref{lem:drift_poly} except that $B_\gen(\step, S)$ depends on both $\step$ and $S$ unlike $\step B_\poly$ in the locally polyhedral case.

Let constants $U_\gen(\step, S), r_\gen(S), \rho_\gen(S), D_\gen(\step, S)$ be
\begin{align}
  U_\gen(\step, S) &\defeq \frac{2 \prts{ D_\gen(\step, S) +1 }}{{r_\gen(S)}^2}, \label{eq:U_gen}\\
  r_\gen(S) &\defeq \frac{3 \Lgen(S)}{24 G^2 + 2 \Lgen(S) G}, \label{eq:r_gen}\\
  \rho_\gen(S) &\defeq 1 - \frac{3 \Lgen(S)^2}{4\prts{24 G^2 + 2 \Lgen(S) G}}, \label{eq:rho_gen}\\
  D_\gen(\step, S) &\defeq \frac{ \prts{ e^{2 G r_\gen(S)} - \rho_\gen(S) } e^{\frac{r_\gen(S)B_\gen(\step,S)}{\step} } }{1 - \rho_\gen(S)}. \label{eq:D_gen}
\end{align}

Define the transient time for a general convex function as
\begin{equation}
  \label{eq:transient_gen}
  T_\gen(\step, S) \defeq \left\lceil \frac{ r_\gen(S) \norm{ \vw_0 - \vw_0^\ast } }{ \step \log(1/ \rho_\gen(S) ) } \right\rceil,
\end{equation}
where $r_\gen(S)$ and $\rho_\gen(S)$ are defined in \eqref{eq:r_gen} and \eqref{eq:rho_gen}.

\begin{lemma}
  \label{lem:expect_gen}
Suppose Assumptions \ref{ass:solution} and \ref{ass:G} hold.  Function $F$ is convex with $\Lgen(S)$ defined in \eqref{eq:LS} for all $S \in (0, S_\text{max})$.  When $t \geq T_\gen(\step, S)$, the following holds
\begin{equation*}
  \expect{ \norm{ \vw_t - \vw_t^\ast }^2 } \leq \step^2 U_\gen(\step, S),
\end{equation*}
where constant $U_\gen(\step, S)$ is defined in \eqref{eq:U_gen}.
\begin{proof}
From Lemma \ref{lem:concentration}, the constants are $\beta = \step\Lgen(S)/2$ and $\gamma = B_\gen(\step, S)$, where \eqref{eq:cond_concentration} holds due to Lemma \ref{lem:drift_gen}.  When $t_0 = 0$, we have $k_{t_0} = K_0 = \norm{ \vw_0 - \vw_0^\ast }$ (with probability 1).  From \eqref{eq:concentration}, it holds for all $t \in \prtc{0, 1, 2, \dotsc}$ that
\begin{align}
  \expect{ e^{\frac{r_\gen(S) K_t}{\step}} } &\leq D_\gen(\step,S) \notag\\
  & \quad+ \prts{ e^{\frac{r_\gen(S) k_{t_0}}{\step} } - D_\gen(\step,S) } \rho_\gen(S)^\prtr{ t - 0 } \notag\\
  &\leq D_\gen(\step,S) + e^{\frac{r_\gen(S) K_0}{\step}}\rho_\gen(S)^t, \label{lem-bound-2-gen}
\end{align}
where constants $r_\gen(S), \rho_\gen(S), D_\gen(\step,S)$ are defined in \eqref{eq:r_gen}--\eqref{eq:D_gen} respectively.  We then show that 
\begin{equation}
  \label{lem-bound-1-gen}
  e^{\frac{r_\gen(S) K_0}{\step}}\rho_\gen(S)^t \leq 1 \quad\quad \forall t \geq T_\gen(\step,S).
\end{equation}
Inequality $e^{\frac{r_\gen(S) K_0}{\step}}\rho_\gen(S)^t \leq 1$ is equivalent to $t \geq \frac{ r_\gen(S) K_0 }{ \step \log(1/ \rho_\gen(S) ) }$.  From the definition of $T_\gen(\step,S)$ in \eqref{eq:transient_gen}, it holds that $T_\gen(\step,S) \geq \frac{ r_\gen(S) K_0 }{ \step \log(1/ \rho_\gen(S) ) }$, and the results \eqref{lem-bound-1-gen} follows.

From \eqref{lem-bound-1-gen}, inequality \eqref{lem-bound-2-gen} becomes
\begin{equation*}
  \expect{ e^{\frac{r_\gen(S) K_t}{\step}} } \leq D_\gen(\step,S) + 1 \quad\quad \forall t \geq T_\gen(\step,S).
\end{equation*}
The Chernoff bounds implies for any $m > 0$ that
\begin{align*}
  \prob{ K_t \geq m } \leq e^{-\frac{r_\gen(S) m }{\step}} \prts{ D_\gen(\step,S) + 1 } ~~ \forall t \geq T_\gen(\step,S).
\end{align*}
Using $\expect{ K_t^2 } = 2 \int_0^\infty m \prob{ K_t \geq m } dm$ and the above bound proves the lemma.
\end{proof}
\end{lemma}

The definition of $T_\gen(\step,S)$ in \eqref{eq:transient_gen} implies that the transient time for a general convex function depends on a step size and the curvature near the unique minimum.  Then, Lemma \ref{lem:expect_gen} implies that $\expect{ \norm{ \vw_t - \vw_t^\ast }^2 } \leq \step^2 U_\gen(\step,S)$ every round $t$ after the transient time.

\subsection{Convergence Rate}
We are now ready to prove the convergence rate of the staggered time averages in Algorithm \ref{alg:main} under a general convex function.

\begin{theorem}
  \label{thm:convergence_gen}
Suppose Assumptions \ref{ass:solution} and \ref{ass:G} hold.  Function $F$ is convex with $\Lgen(S)$ defined in \eqref{eq:LS} for all $S \in (0,S_\text{max})$.  It holds for any $t \geq T_\gen(\step, S)$ and any positive integer $T$ that
\begin{equation}
  \label{eq:F_gen}
  \expect{ F(\bar{\vw}_{t}^T) } - F^\ast \leq \frac{\step G^2}{2} + \frac{\step U_\gen(\step,S)}{2T},
\end{equation}
where constant $U_\gen(\step, S)$ is defined in \eqref{eq:U_gen}
\begin{proof}
The theorem follows from \eqref{eq:avgbasic} and Lemma \ref{lem:expect_gen} that
\begin{multline*}
  \expect{ F(\bar{\vw}_{t}^T) } - F^\ast \leq \frac{\step G^2}{2} + \frac{\step^2 U_\gen(\step,S)}{2\step T}\\ \forall t \geq T_\gen(\step, S).
\end{multline*}
\end{proof}
\end{theorem}

The transient time \eqref{eq:transient_gen} and Theorem \ref{thm:convergence_gen} can be interpreted as a class of convergence bounds that can be optimized over any $S \in (0, S_\text{max})$.  Indeed, the values of $S$ near the minimum of $F$ plays a crucial role in \eqref{eq:B_gen}, which affects much of the analysis in this section.

\section{Fast convergence for deterministic problems} 
\label{sec:deterministic}
This section revisits problems with the locally polyhedral structure, so that Assumptions \ref{ass:solution}, \ref{ass:G}, \ref{ass:poly} hold.  However, it considers a \emph{deterministic} scenario where the oracle returns the exact subgradient, rather than an unbiased stochastic subgradient.   It is shown that a variation on the basic algorithm that uses a variable step size can significantly improve the convergence rate.   Specifically, fix $\epsilon > 0$.  The basic algorithm of Section \ref{sec:poly} produces an $O(\epsilon)$-approximation within $O(1/\epsilon)$ rounds.  The modified algorithm of this section does the same $O(\epsilon)$-approximation with only $O( \log(1/\epsilon))$ rounds.  In particular, this is faster than the lower bound $\Omega(1/\sqrt{T})$ for a non-smooth function with Lipschitz continuity in \cite{Nesterov:Intro}. This does not contradict the Nesterov result in \cite{Nesterov:Intro}, because that result shows the existence of a function with $\Omega(1/\sqrt{T})$ convergence rate, while the locally polyhedral structure does not fall into a class of that function.  Interestingly, the algorithm in this section is Faster than other algorithms with $O(1/T^2)$ convergence rates \cite{Nesterov:Intro,Tseng:Fast_gradient}.

Assume the function $F(\vw)$ is Lipschitz continuous over $\vw \in \set{W}$ with Lipschitz constant $H > 0$, so that:
\begin{equation}
  \label{ass:Lipschitz}
  \abs{ F(\vw) - F(\vw') } \leq H \norm{ \vw - \vw' } \quad \forall \vw, \vw' \in \set{W}.
\end{equation}
Assume there is a known positive value $ Z < \infty$ such that $\norm{\vw_0 - \vw_0^\ast} \leq Z$.  Fix $\epsilon > 0$, and fix $M$ as any positive integer.  The idea is to run the algorithm over successive frames.  Label the frames $i \in \prtc{1, 2, 3, \dotsc, M}$.  Let $\vw_{[i]}$ be the initial vector in $\set{W}$ at the start of frame $i$.   Define $\vw_{[1]} = \vw_0$.  Define constants $U_\poly$, $r_\poly$, $\rho_\poly$, $D_\poly$ as in \eqref{eq:U_poly}--\eqref{eq:D_poly}.  Define $\theta = \max\prts{\sqrt{U_\poly}, Z}$, and define the \emph{frame size} $T$ as: 
\begin{equation*}
  T = \left\lceil \frac{ 2 r_\poly \theta }{ \log(1/\rho_\poly)} \right\rceil
\end{equation*}
The algorithm for each frame $i \in \prtc{1, 2, 3, \dotsc}$ is: 
\begin{itemize} 
\item Define the step size for frame $i$ as $\alpha_{[i]} = 2^{-i}$. 
\item Run the constant step size algorithm \eqref{alg:subgrad} using step size $\alpha_{[i]}$ over $T$ rounds, using initial vector $\vw_{[i]}$.
\item Define $\vw_{[i+1]}$ as the $\vw_t$ vector computed in the last round of frame $i$.
\end{itemize} 

Notice that the completion of $M$ frames requires $MT = O(M)$ rounds. The vector computed in the last round of the last frame is defined as $\vw_{[M+1]}$.  The next theorem shows that this vector is indeed an $O(2^{-M})$-approximation.

\begin{theorem} 
In the deterministic setting and when Assumptions \ref{ass:solution}, \ref{ass:G}, \ref{ass:poly} hold, the final vector $\vw_{[M+1]}$ satisfies: 
\begin{align}
\norm{\vw_{[M+1]} - \vw_{[M+1]}^\ast} &\leq \theta 2^{-M} \label{eq:end-distance}\\
F(\vw_{[M+1]}) - F^\ast &\leq  \theta H 2^{-M} \label{eq:end-fval}
\end{align}

\begin{proof}  The proof is by induction on the rounds $i \in \prtc{1, 2, \dotsc, M}$.  Assume the following holds on a given $i \in \prtc{1, 2, \dotsc, M}$: 
\begin{equation} 
  \label{eq:norm-bound-frame} 
  \norm{\vw_{[i]} - \vw_{[i]}^\ast} \leq \theta 2^{-(i-1)}
\end{equation} 
This holds by assumption on the first frame $i=1$. We now show \eqref{eq:norm-bound-frame} holds for $i+1$.  The goal is to use Lemma \ref{lem:expect_poly} with initial condition $\vw_{[i]}$.  Since the step size is $\alpha_{[i]}$ for this frame, the value $T_\poly(\alpha_{[i]})$ defined in \eqref{eq:T_poly} satisfies: 
\begin{align*}
  T_P(\alpha_{[i]}) = \left\lceil \frac{r_\poly\norm{\vw_{[i]}-\vw_{[i]}^\ast}}{\alpha_{[i]}\log(1/\rho_\poly)}\right\rceil \leq \left\lceil \frac{r_\poly\theta 2^{-(i-1)}}{\alpha_{[i]}\log(1/\rho_\poly)}\right\rceil = T
\end{align*}
where the inequality holds by the induction assumption \eqref{eq:norm-bound-frame}, and the last equality holds by definition of $\alpha_{[i]}$.   Recall that $\vw_{[i+1]}$ is defined as the final $\vw_t$ value after $T$ rounds of the frame. It follows by Lemma \ref{lem:expect_poly} that: 
\begin{equation*}
  \expect{\norm{\vw_{[i+1]} - \vw_{[i+1]}^\ast}^2} \leq \alpha_{[i]}^2 U_\poly \leq \alpha_{[i]}^2 \theta^2
\end{equation*}
On the other hand, this deterministic setting produces a deterministic sequence, so that all expectations can be removed: 
\begin{equation*}
  \norm{\vw_{[i+1]}- \vw_{[i+1]}^\ast}^2 \leq \alpha_{[i]}^2 \theta^2
\end{equation*}
Taking a square root and using the definition of $\alpha_{[i]}$ proves: 
\begin{equation*}
 \norm{\vw_{[i+1]}-\vw_{[i+1]}^\ast} \leq \theta 2^{-i}
\end{equation*}
This completes the induction, so that \eqref{eq:norm-bound-frame} holds for all $i \in \prtc{1, 2, \dotsc, M+1}$. Substituting $i=M+1$ into \eqref{eq:norm-bound-frame} proves \eqref{eq:end-distance}.  The inequality \eqref{eq:end-fval} follows from \eqref{eq:end-distance} and the Lipschitz property \eqref{ass:Lipschitz}.
\end{proof} 
\end{theorem}

\section{Experiments}
\label{sec:sim}
In this section, Algorithm \ref{alg:main} (``Staggered'') is compared to the polynomial-decay averaging (``Polynomial'') in \cite{Ohad:averaging}.  We also proposed another heuristic algorithm (``Heuristic''), which has a promising convergence rate.  This heuristic algorithm replaces $\step$ in Algorithm \ref{alg:main} with 
\begin{equation*}
  \step_t = \max\prts{ \step, c/(t+1) } \quad\quad \forall t \in \prtc{0, 1, 2, \dotsc}
\end{equation*}
where $c$ is some real-valued positive constant.  This modification does not change the convergence rates in Sections \ref{sec:poly} and \ref{sec:gen}, because it only adds $O(1/\step)$ rounds into the previous bounds.

For the purpose of comparison, Algorithm \ref{alg:main} uses $\step = 10^{-4}$.  However, higher accuracy can be achieved by a smaller step size.  The heuristic algorithm sets $c = 1$.  The polynomial-decay averaging algorithm uses $c = 1$ and $\eta = 3$ (defined in \cite{Ohad:averaging}).  Note that the stochastic subgradient algorithm with a constant step size (``Constant'') is the by product of Algorithm \ref{alg:main}.

\begin{figure}
  \centering
  \includegraphics[scale=0.5]{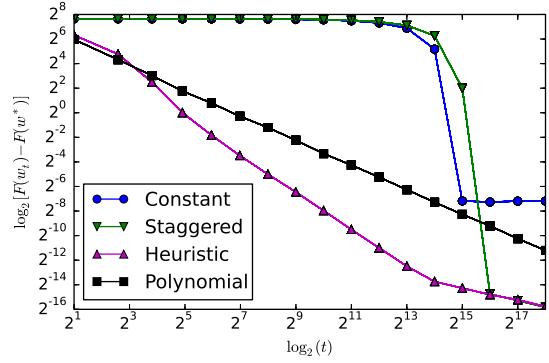}
  \caption{Results of algorithms and a locally polyhedral function}
  \label{fig:LP_unique}
\end{figure}

A locally polyhedral function $F = \norm{ \vw }_1$ is considered where $\vw \in [-4, 4]^{100}$.  When $\vg(\vw)$ is a subgradient of $F$ at $\vw$, a stochastic subgradient is $\hat{\vg}(\vw) = \vg(\vw) \times X$ where $X$ is a uniform random variable from $0$ to $2$, so $\expect{\vg(\vw)|\vw} = \vg(\vw)$.  Ten experiments are performed, and the average values are sampled at $\prtc{2^k-2}_{k=1}^{18}$ (one round before Algorithm \ref{alg:main} resets the averages).  Results are shown in Figure \ref{fig:LP_unique}.  Both axes of Figure \ref{fig:LP_unique} are in a log scale.

The plots of Algorithm \ref{alg:main} and the polynomial-decay algorithm cross each other, because the former has faster convergence rate.  The subgradient algorithm with constant step size stops improving due to the fixed value of the step size.  However, Algorithm \ref{alg:main} keeps improving after the stop.  This can be explained by \eqref{eq:F_poly} where the average helps reducing the last term on the right-hand-side.  The plot of the heuristic algorithm shows its convergence.

\begin{figure}
  \centering
  \includegraphics[scale=0.5]{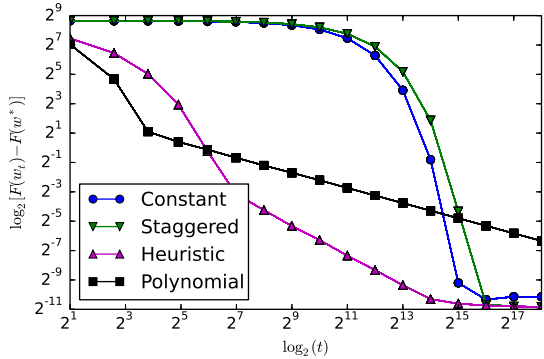}
  \caption{Results of algorithms and a general convex function}
  \label{fig:GP_unique}
\end{figure}

A general convex function $F(\vw) = \sum_{i=1}^{100} F_i(w^{(i)})$ is considered where, for $i \in \prtc{1, \dotsc, 100}$, $w^{(i)}$ is the $i$-component of $\vw$ and
\begin{equation*}
  F_i(w^{(i)}) = \left\{
    \begin{array}{rl}
      -w^{(i)}, & \text{if}~ w^{(i)} < 0 \\
      (w^{(i)})^2, & \text{if}~ w^{(i)} \geq 0.
    \end{array} \right.
\end{equation*}
The $i$-component of a stochastic subgradient is $\hat{\vg}_i(\vw) = \vg_i(\vw) + Y$, where $\vg_i(\vw)$ is the $i$-component of the true subgradient of $F$ at $\vw$ and $Y$ is a uniform random variable between -1 and 1.  Simulation uses the same parameters as the locally polyhedral case.  The results in Figure \ref{fig:GP_unique} have a similar trend as in the locally polyhedral case except that the plot of the stochastic subgradient algorithm with constant step size crosses the plot of the polynomial decay averaging.

\begin{figure}
  \centering
  \includegraphics[scale=0.5]{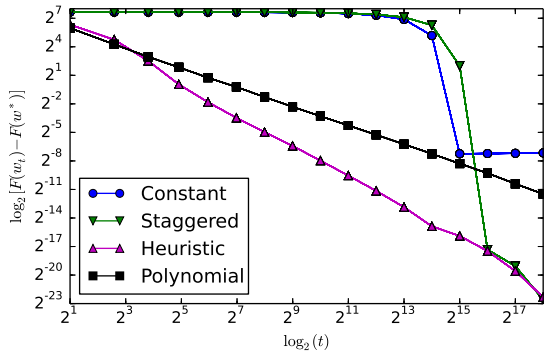}
  \caption{Results of algorithms and the function \eqref{eq:F_LP_non_unique}}
  \label{fig:LP_non_unique}
\end{figure}

Then we consider a non-smooth convex function $F(\vw) = \sum_{i=1}^{100} F_i(w^{(i)})$ where
\begin{equation}
  \label{eq:F_LP_non_unique}
  F_i(w^{(i)}) = \left\{
    \begin{array}{rl}
      w^{(i)} - \frac{10^{-6}}{2},  & \text{if}~ w^{(i)} \geq \frac{10^{-6}}{2} \\
      -w^{(i)} - \frac{10^{-6}}{2}, & \text{if}~ w^{(i)} \leq -\frac{10^{-6}}{2} \\
      0,                  & \text{otherwise}.
    \end{array} \right.
\end{equation}
This function has uncountable minimizers.  The stochastic subgradient is the component-wise addition of the true subgradient and the uniform random variable $Y$.  Simulation results are shown in Figure \ref{fig:LP_non_unique}.  Comparing these results to the results in Figure \ref{fig:LP_unique} shows the same trend of  convergence rates even though function $F$ in \eqref{eq:F_LP_non_unique} does not satisfy the uniqueness assumption.

\section{Conclusion}
\label{sec:conclusion}
This paper considers stochastic non-smooth convex optimization.  We propose the staggered time average algorithm and prove its performance.  When a function with a unique minimum satisfies the locally polyhedral structure, the algorithm has $O(1/T)$ convergence rate.  For a general convex function with a unique minimum, we derive a class of bounds on the convergence rate of the algorithm.  For a special case of deterministic problems with the locally polyhedral structure, an algorithm with $O(1/\epsilon^{1/M})$ convergence is proposed.

\bibliographystyle{IEEEtran}
\bibliography{Reference}

\end{document}

%% file: SuchaLib.tex
\newcommand{\set}[1]{ {\mathcal{#1}} }
\newcommand{\vect}[1]{ {\boldsymbol{#1}} }
\newcommand{\expect}[1]{ {\mathbb{E}\left[ {#1} \right]} }
\newcommand{\prob}[1]{ {\mathbb{P}\left\{ {#1} \right\}} }

\newcommand{\abs}[1]{ {\left| {#1} \right|} }

\newcommand{\norm}[1]{ {\left\lVert {#1} \right\rVert} }

\newcommand{\inner}[2]{ {\left\langle {#1} , {#2} \right\rangle} }

\newcommand{\RealSet}{ {\mathbb{R}} }

\newcommand{\defeq}{ {\triangleq} }

\DeclareMathOperator{\arginf}{arginf}
\DeclareMathOperator{\Arginf}{Arginf}

\newcommand{\prts}[1]{ {\left[ {#1} \right]} }
\newcommand{\prtr}[1]{ {\left( {#1} \right)} }
\newcommand{\prtc}[1]{ {\left\{ {#1} \right\}} }